\documentclass[12pt,leqno]{article}

\usepackage{amssymb,amsmath,amsthm}
\usepackage{mathrsfs}
\numberwithin{equation}{section}

\newcommand{\corps}{\operatorname{\mathbf{k}}}

\newcommand{\Int}{{\rm Int}}

\newcommand{\supp}{{\rm supp}}

\newcommand{\lan}{\langle}
\newcommand{\ran}{\rangle}
\newcommand{\bra}[1]{\lan{#1}\ran}
\newcommand{\pset}[2]{\{#1\,;\,\mbox{#2}\}}
\newcommand{\eps}{\varepsilon}
\newcommand{\Supp}{\operatorname{Supp}}
\newcommand{\dist}{{\rm dist}}
\newcommand{\Ss}{\operatorname{SS}}
\newcommand{\codim}{\operatorname{codim}}

\def\etens{\mathbin{\boxtimes}}
\newcount\temp
\temp=\catcode`\@
\catcode`\@=11
\def\@bletens{\mathbin{\etens^{L}}}
\def\@letens_#1{\mathbin{\etens_{\raise1.5ex\hbox to-.1em{}#1}^{L}}}
\def\letens{\@ifnextchar _{\@letens}{\@bletens}}
\catcode`\@=\temp

\def\phi{{\varphi}}
\def\epsilon{\varepsilon}

\newcommand{\md[1]}{{\rm Mod}({#1})}

\newcommand{\ba}{\begin{array}}
\newcommand{\ea}{\end{array}}

\newcommand{\C}{\mathbb{C}}

\newcommand{\R}{\mathbb{R}}

\newcommand{\Z}{\mathbb{Z}}

\renewcommand{\to}[1][\rule{.7em}{0em}]{\xrightarrow{#1}}

\newcommand{\isoto}[1][]{\xrightarrow[#1]{\sim}}

\renewcommand{\hom}[1][]{{\mathcal{H}\kern-1pt{om}_{\raise1.5ex\hbox to.1em{}#1}}}

\newcommand{\rhom}[1][]{{R{\mathcal{H}}{om}_{\raise1.5ex\hbox to.1em{}#1}}}
\newcommand{\ext}[1][]{{\mathcal{E}xt}_{\raise1.5ex\hbox to.1em{}#1}}

\newcommand{\Hom}[1][]{\mathrm{Hom}_{\raise1.5ex\hbox to.1em{}#1}}
\newcommand{\RHom}[1][]{\mathrm{RHom}_{\raise1.5ex\hbox to.1em{}#1}}

\theoremstyle{plain}

\newtheorem{theorem}{Theorem}[section]
\newtheorem{proposition}[theorem]{Proposition}

\newtheorem{corollary}[theorem]{Corollary}

\theoremstyle{definition}

\newtheorem{definition}[theorem]{Definition}

\newtheorem{example}[theorem]{Example}

\newtheorem{remark}[theorem]{Remark}

\newcommand{\eq}{\begin{eqnarray}}
\newcommand{\eneq}{\end{eqnarray}}
\newcommand{\eqn}{\begin{eqnarray*}}
\newcommand{\eneqn}{\end{eqnarray*}}

\newenvironment{nnum}{
  \begin{enumerate}
  \itemsep=0pt
  
  }
  {\end{enumerate}}

\newenvironment{anum}{
  \begin{enumerate}
  \itemsep=0pt
  
  }
  {\end{enumerate}}

\newcommand{\bnum}{\begin{nnum}}
\newcommand{\enum}{\end{nnum}}
\newcommand{\banum}{\begin{anum}}
\newcommand{\eanum}{\end{anum}}
\newcommand{\ol}{\overline}
\def\dddt{{\raise-.3em\hbox{$\big\cdot$}}}

\newcommand{\cl}{\colon}
\newcommand{\hs}[1]{\hspace*{#1}}

\newcommand{\bl}{\bigl}
\newcommand{\br}{\bigr}

\newcommand{\Rc}{{\mathbb{R}-c}}
\newcommand{\Cc}{{\mathbb{C}-c}}

\begin{document}

\author{Masaki Kashiwara,
Teresa Monteiro Fernandes%
\\ 
and Pierre Schapira}
\title{Involutivity of truncated microsupports}

\maketitle
\footnote{
The research of the first author is partially supported by
Grant-in-Aid for Scientific Research (B) 13440006,
Japan Society for the Promotion of Science.}
\footnote{
The research of the second author
was supported by FCT and Programa Ci{\^e}ncia,
Tecnologia e Inova\c c{\~a}o do Quadro Comunit{\'a}rio de Apoio.}

\begin{abstract}
Using a result of J-M. Bony, we prove
the weak involutivity of truncated microsupports. More precisely,
given a  sheaf $F$ on  a real manifold and $k\in\Z$, if
 two functions vanish on $\Ss_k(F)$, then so
does their Poisson bracket. 
\end{abstract}

\section{Introduction}\label{section:intro}
The notion of microsupport of sheaves
was introduced by two of the present authors (M.K.\ and P.S.) 
in the course of the study 
of the theory of linear partial differential equations. 
References are made to \cite{K-S1}.
These
authors also introduced a variant of this notion, that of ``truncated 
 microsupport'' and developed its study with the third author
 in \cite{KMS}.

A crucial result in the microlocal theory of sheaves
is the involutivity of the microsupport. This property 
does not hold true any more for the truncated microsupport, but a weak
form of it does. More precisely, 
we prove here that
the truncated microsupport is stable by Poisson bracket, that is, 
if two functions vanish on it, so
does their Poisson bracket. 
The main technical tool is a result of J.-M. Bony which asserts the
same property for the normal cone to a closed subset.

\section{Notations and review}\label{S:1}
We will follow the notations of \cite{K-S1} and \cite{KMS}.
For the reader convenience, we recall some of them as well as 
the definition of the truncated
microsupport.
Let $X$ be a real analytic manifold. We denote by
$\tau\colon TX\to X$ the tangent bundle to $X$ and
by $\pi\colon T^*X\to X$ the
 cotangent bundle. 
For a smooth submanifold $Y$ of $X$, $T_YX$ denotes
 the normal bundle to $Y$ and $T^*_YX$ the conormal bundle.
In particular, $T^*_X X$ is the zero section of $T^*X$.
We set
$\dot{T}^*X=T^*X\setminus X$, and denote by
$\dot{\pi}\colon
\dot{T}^*X\to X$  the restriction of $\pi$ to $\dot{T}^*X$.

For a morphism $f\colon X\to Y$ of real manifolds, we denote
by
$$f_{\pi}\colon X\times_YT^*Y\to T^*Y
\mbox{ and }f_d \colon X\times_YT^*Y\to T^*X$$
the associated morphisms.

For a subset $A$ of $T^*X$,
we denote by $A^a$ the image of $A$ by the
antipodal map $a\colon (x;\xi)\mapsto (x;-\xi)$.
The closure of $A$ is denoted by  $\overline{A}$.

For a cone $\gamma\subset TX$,
the polar cone $\gamma^\circ$ to $\gamma$ is the convex cone in $T^*X$
defined by
$$\gamma^\circ=\{(x;\xi)\in T^*X;
\mbox{$x\in\pi(\gamma)$ and $\langle v,\xi\rangle\geq 0$ for any
$(x;v)\in \gamma$}\}.$$
A closed convex cone
$\gamma$ is called {\em proper}
if $0\in\gamma$ and $\Int(\gamma^\circ)\neq\emptyset$.

Let $\corps$ be a field.
One denotes by $\md[\corps_X]$ the abelian category of sheaves of 
$\corps$-vector spaces and by $D^b(\corps_X)$ its bounded derived
category. One denotes by $D^b_\Rc(\corps_X)$ the full triangulated subcategory
of $D^b(\corps_X)$ consisting of objects with $\R$-constructible cohomology.
If $X$ is a complex manifold, one denotes by 
$D^b_\Cc(\corps_X)$ the full triangulated subcategory
of $D^b(\corps_X)$ consisting of objects with $\C$-constructible cohomology.

If $S$ is a locally closed subset of $X$, one denotes by
$\corps_{XS}$ the sheaf on $X$ which is the
constant sheaf with stalk $\corps$ on $S$ and $0$ on $X\setminus S$.
If there is no risk of confusion, we may write $\corps_{S}$ instead
of $\corps_{XS}$.

For $k\in \Z$, we denote as usual by  $D^{\geq k}(\corps_X)$
(resp.\ $D^{\leq k}(\corps_X)$)
the full additive subcategory of $D^b(\corps_X)$ consisting of
objects $F$ satisfying $H^j(F)=0$ for any $j<k$
(resp.\ $H^j(F)=0$ for any $j>k$).

We denote by $\tau^{\leq k}\cl D(\corps_X)\to D^{\leq k}(\corps_X)$
the truncation functor. Recall that for $F\in D(\corps_X)$
the morphism $\tau^{\le k}F\to F$
induces isomorphisms
$H^j(\tau^{\le k}F)\isoto H^j(F)$
for $j\le k$ and $H^j(\tau^{\le k}F)=0$
for $j>k$.

If $F$ is an object of $D^b(\corps_X)$,
$\Ss(F)$ denotes its microsupport, a closed
$\R^{+}$-conic subset of $T^*X$. 
For $p\in T^*X$, $D^b(\corps_X; p)$ denotes the localization of 
$D^b(\corps_X)$ by the full triangulated subcategory consisting of objects $F$
such that $p\notin \Ss(F)$. A property holds ``microlocally on a subset $S$
of $T^*X$'' if it holds in the category $D^b(\corps_X; p)$ for any
$p\in S$.

We recall the definition of involutivity of \cite{K-S1}. This notion 
  makes use of that of
``normal cone'' (\cite[Definition 4.1.1]{K-S1}).
For a pair of subsets $S_1$ and $S_2$ of a manifold $X$,
the normal cone $C_p(S_1,S_2)$ at $p\in X$ is defined as follows:
it is a closed cone in the tangent space $T_pX$
consisting of points $v$ such that,
for a local coordinate system, there exist a sequence $\{x_n\}_{n}$ in $S_1$,
$\{y_n\}_n$ in $S_2$ and a sequence $\{a_n\}_n$ in $\R_{\ge0}$
such that $x_n$ and $y_n$ converge to $p$
and $a_n(x_n-y_n)$ converges to $v$.
For a subset $S$, $C_p(S,\{p\})$ is denoted by $C_p(S)$.

\begin{definition}\label{def:invol}(\cite[Definition 6.5.1]{K-S1})
Let $S$ be a locally closed subset of $T^*X$ and let $p\in S$. One
says that $S$ is involutive at $p$ if for any $\theta\in T_p^*(T^*X)$
such that the normal cone $C_p(S,S)$ is contained in the hyperplane
$\pset{v\in T_p(T^*X)}{$\langle v,\theta\rangle=0$}$ one has: 
$H(\theta)\in C_p(S)$. Here $H\cl T^*_p(T^*X)\isoto T_p(T^*X)$
is the Hamiltonian isomorphism.

If $S$ is involutive at each $p\in S$, one says that $S$ is involutive.
\end{definition}
The involutivity theorem of \cite[Theorem 6.5.4]{K-S1} asserts that the
microsupport of sheaves is involutive.

The following definition was introduced by the authors
of \cite{K-S1} and developed in \cite{KMS}. 

\begin{definition}\label{def:ssk}
Let $X$ be a real analytic manifold and let $p\in T^*X$.
Let $F\in D^b(\corps_X)$ and $k\in\Z$.
The closed conic subset $\Ss_k(F)$ of $T^*X$ is defined by:
$p\notin \Ss_k(F)$ if and only if the following condition is satisfied.
\bnum
\item[{\rm (i)}]
There exists an open conic neighborhood $U$
of $p$ such that for any
$x\in \pi(U)$ and for any $\R$-valued C${}^1$-function $\phi$
defined on a neighborhood of $x$ such
that $\phi(x)=0$, $d\phi(x)\in U$, one has
\begin{equation}\label{E:1}
H^j_{\{\phi\geq 0\}}(F)_x=0\quad\mbox{for any $j\le k$.}
\end{equation}
\enum
\end{definition}
We refer to \cite{KMS} for equivalent definitions. In particular, it
is proved in loc.\ cit.\ that one can replace 
the condition that $\phi$ is of class  ${\rm C}^1$ by $\phi$ is of 
class ${\rm C}^\alpha$, where 
$\alpha\in\Z_{\ge1}\cup\{\infty,\omega\}$.

This condition is also equivalent to:
\bnum
\item[{\rm (ii)}]
there exists $F'\in D^{\geq k+1}(\corps_X)$
and an isomorphism $F\simeq F'$ in $D^b(\corps_X;p)$.
\enum

\section{Normal cone}\label{S:2}

Let us recall the notion of
the normal vector due to J.-M.~Bony.
Here we set for $r\ge0$
\eqn
&&B_r(x)
=\pset{y\in\R^n}{$\Vert y-x\Vert<r$},
\eneqn
the open ball with center $x$ and radius $r$.

\begin{proposition}\label{pro:normal}
Let $X$ be a real manifold,
$S$ a closed subset of $X$, and $\alpha\in\Z_{\ge1}\cup\{\infty,\omega\}$.
Then 
the following subsets of $T^*X$ are equal.
\bnum
\item[{\rm (i)}] $\Ss_0(\corps_S)$.
\item[{\rm (ii)}${}_\alpha$]
The closure of the set of points $p=(x;\xi)\in T^*X$
such that $x\in S$ and there is a
{\rm C}${}^{\alpha}$ function $\phi$ defined on a neighborhood $U$ of $x$ 
such that $\phi(x)=0$, $d\phi(x)=\xi$ 
and $S\cap U\subset\phi^{-1}(\R_{\ge0})$.
\item[{\rm (iii)}]
The closure of the set of points $p=(x;\xi)\in T^*X$
such that $x\in S$ and $C_x(S)\subset\{v\in T_{x}X;\bra{v,\xi}\ge0\}$.
\enum
If $X$ is an open subset of $\R^n$, then
the above sets are equal to the following set.
\bnum
\item[{\rm (iv)}]
The closure of the set of points $p=(x;\xi)\in T^*X$
such that $x\in S$ and  
the open ball $B_{t\Vert\xi\Vert}(x-t\xi)$
does not intersect $S$ for some $t>0$.
\enum
\end{proposition}

\begin{proof}
We may assume that $X=\R^n$.

\noindent
(i)=(ii)${}_\alpha$ is an immediate consequence of Definition \ref{def:ssk}.
In particular, the set (ii)${}_\alpha$ does not depend on $\alpha$.

\noindent
(ii)${}_\alpha$=(iv) for $\alpha\ge2$ is obvious.

\noindent
(ii)${}_\alpha$ $\subset$ (iii) is clear.

\noindent
(iii) $\subset$ (iv).
Let $p=(x_0;\xi_0)\in T^*X$ with
$x_0\in S$ and $C_{x_0}(S)\subset\{v\in T_{x}X;\bra{v,\xi_0}\ge0\}$.
We have to prove that $p$ belongs to the set (iv).
The proof is very similar to that of Lemma 3.3 of \cite{KMS}, but
for the sake of completeness, we shall repeat it.

If $\xi_0=0$, then it is trivially
true. Hence we may assume that $\xi_0\not=0$.

We shall show that for an arbitrary open conic neighborhood
$U$ of $p$, there exists a point $(x_1;\xi_1)\in U$
such that
the open ball $B_{\Vert\xi_1\Vert}(x_1-\xi_1)$
satisfies
\eq\label{eq:closure}
&&x_1\notin\ol{B_{\Vert\xi_1\Vert}(x_1-\xi_1)\cap S}.
\eneq
This implies $B_{t\Vert\xi_1\Vert}(x_1-t\xi_1)\cap S=\emptyset$
for $0<t\ll1$.

Let us take an open neighborhood
$V$ and a proper closed convex cone $\gamma$
such that $\xi_0\in \Int(\gamma^\circ)$,
$\ol{V}\times (\gamma^\circ\setminus\{0\})\subset U$ and
$\Int(\gamma)\not=\emptyset$.
Since $C_{x_0}(S)\cap\gamma^a\subset\{0\}$,
there is $\rho>0$ such that
$H_-:=\pset{x\in\R^n}{$\bra{x-x_0,\xi_0}>-\rho$}$
satisfies $\ol{H_-}\cap(x_0+\gamma^a)\subset V$ and 
$S\cap \ol{H_-}\cap(x_0+\gamma^a)=\{x_0\}$.
Set $S_0:=S\cap\ol{H_-}$.

Let us define the function $\psi(x)$ on $\R^n$ by
$\psi(x)=\dist(x,\gamma^a):=\inf\{\Vert y-x\Vert\,;\,y\in\gamma^a\}$.
It is well known that
$\psi$ is a continuous function on $\R^n$,
and $C^1$ on $\R^n\setminus\gamma^a$.
More precisely for any $x\in \R^n\setminus\gamma^a$,
there exists a unique
$y\in\gamma^a$ such that
$\psi(x)=\Vert x-y\Vert$.
Moreover $d\psi(x)=\Vert x-y\Vert^{-1}(x-y)\in\gamma^\circ{}\setminus\{0\}$.
Furthermore $B_{\psi(x)}(y)$
is contained in $\{z\in\R^n;\psi(z)<\psi(x)\}$.

For $\epsilon>0$, we set $\gamma^a_\epsilon=\{x\in\R^n;\psi(x)<\epsilon\}$.
Then $\gamma^a_\epsilon$ is an open convex set.
Moreover $\gamma^a_\eps+\gamma^a=\gamma^a_\eps$.

Let us take $v\in\Int(\gamma)$,
and $\delta>0$ such that
$x_0-\delta v\in H_-$.
Then take $\eps>0$ such that
$\gamma^a_\eps-\delta v\subset\Int(\gamma^a)$,
$(x_0+\ol{\gamma^a_\eps})\cap\ol{H_-}\subset V$
and $(x_0+\ol{\gamma^a_\eps})\cap S_0\subset H_-$.

Set $W_t=x_0+\gamma_{\eps}^a+tv$ for $t\in\R$.
Then one has
\eq
&&\mbox{$W_t=\bigcup_{t'<t}W_{t'}$ and 
$\ol{W_t}=\bigcap_{t'>t}{W_{t'}}=\bigcap_{t'>t}\ol{W_{t'}}$,}
\label{eq:clos}\\
&&\mbox{$x_0\in W_t\cap S$ for $t\ge0$,
and $W_t\cap \ol{H_-}=\emptyset$ for $t\ll0$.}\\
&&\mbox{$\ol{W_t}\cap S_0\subset H_-\cap V$ for $t\le0$}
\eneq
 Hence, for any closed subset $K$ of $\ol{H_-}$
and $t\in\R$ such that $K\cap\ol{W_t}=\emptyset$,
there exists $t'>t$ such that
$K\cap{W_{t'}}=\emptyset$.

Let us set $c=\sup\{t;W_t\cap S_0=\emptyset\}$.
Then $c\le 0$ and
$W_c\cap S_0=\emptyset$. 
By the above remark, one has 
$\ol{W_c}\cap S_0\not=\emptyset$.
Take $x_1\in S_0\cap\partial W_c$.
Here $\partial W_c:=\ol{W_c}\setminus W_c$ is the boundary of
$W_c$.
The point $x_1$ belongs to $H_-\cap V$.

As seen before, there exists
an open ball
$B_{\eps}(y)$ such that
$B_{\eps}(y)\subset W_c$, $\Vert x_1-y\Vert=\eps$
and $x_1-y\in\gamma^\circ$.
Hence $B_\eps(y)\cap S\cap H_-=\emptyset$, 
while $H_-$ is a neighborhood of $x_1$.
Hence $(x_1;x_1-y)$ belongs to $U$ and  satisfies condition 
\eqref{eq:closure}.
\end{proof}

\begin{definition}
For a closed subset $S$ of $X$, the closed subset
given in Proposition \ref{pro:normal} is called the {\em $0$-conormal cone}
of $S$ and denoted by $N^*_0(S)$.
\end{definition}

Note that $N^*_0(S)$ is a closed cone and
it satisfies $\pi_X(N^*_0(S))=\pi_X(N^*_0(S)\cap T^*_XX)=S$.
Note also that one has
$N^*_0(S)=\Ss_0(\corps_S)\subset\Ss(\corps_S)$.

\begin{example}
\bnum
\item
If $S=X$, then $N^*_0(S)=T^*_XX$.
\item
If $S$ is a closed submanifold, then $N^*_0(S)=T^*_SX$.
\item
If $X=\R^2$ and $S=\pset{(x,y)\in X}{$x\ge 0$ or $y\ge0$}$.
Then
\eqn
&&N^*_0(S)=\pset{(x,y;0,0)\in T^*_XX}{$(x,y)\in S$}\\
&&\hs{50pt}\cup\pset{(x,y;\xi,\eta)\in T^*X}%
{$y=0$, $x\le0$, $\eta\ge0$ and $\xi=0$}\\
&&\hs{50pt}
\cup\pset{(x,y;\xi,\eta)\in T^*X}%
{$x=0$, $y\le0$, $\xi\ge0$ and $\eta=0$}.
\eneqn
One has the microlocal isomorphisms 
\eqn
\corps_{XS}\simeq
\left\{
\ba{ll}
\corps_X&\mbox{on $\pset{(x,y;0,0)}{$x>0$ or $y>0$}$,}\\
\corps_{\{y=0\}}&\mbox{on $\pset{(x,0;0,\eta)}{$x<0$ and $\eta>0$}$,}\\
\corps_{\{x=0\}}&\mbox{on $\pset{(0,y;\xi,0)}{$y<0$ and $\xi>0$}$,}\\
\corps_{\{(0,0)\}}[-1]&\mbox{on $\pset{(0,0;\xi,\eta)}{$\xi>0$ and $\eta>0$}$.}
\ea\right.
\eneqn
Note that (see Theorem \ref{th:Rcons} below)
$$
\Ss(\corps_{XS})=\Ss_1(\corps_{XS})=N^*_0(S)\cup\pset{(0,0;\xi,\eta)\in T^*X}
{$\xi,\eta\ge0$}
$$
and this set is different from $N^*_0(S)$.

\enum
\end{example}

\begin{proposition}
\bnum
\item
Let $f\cl X\to Y$ be a morphism of manifolds and let 
$S$ be a closed subset of $X$.
Assume that $f(S)$ is a closed subset of $Y$.
Then one has
$N^*_0(f(S))\subset f_\pi f_d^{-1}(N^*_0(S))$.
The equality holds in case $f$ is a closed embedding.
\item
For a closed subset $S_1$ of $X_1$
and a closed subset $S_2$ of $X_2$, one has
$N^*_0(S_1\times S_2)=N^*_0(S_1)\times N^*_0(S_2)$.
\item 
If a closed subset $S$ of $X$ satisfies
$N^*_0(S)\subset T^*_XX$, then $S$ is an open subset.
\item
Let $f$ be a C$\,{}^1$-function on $X$, let $S$ a closed subset,
and let $c\in\R$.
Assume that $f\vert_S\cl S\to\R$ is proper, $f(S)\subset [a,\infty)$ for $a\ll0$.
If $df(x)\notin N^*_0(S)$ for $x$ such that $f(x)<c$,
then $f(S)\subset [c,\infty)$.
\item
Let $\gamma$ be a proper closed convex cone of $\R^n$,
$\Omega$ an open subset of $\R^n$ such that
$\Omega+\gamma^a=\Omega$,
and $S$ a closed subset of $\Omega$ such that
$S$ is relatively compact in $\R^n$.
If $N^*_0(S)\cap (S\times\Int(\gamma^\circ))=\emptyset$,
then $S$ is an empty set.
\enum
\end{proposition}
\begin{proof}
(i), (ii), (iii) are easy exercises.

\medskip
(iv) If $f\vert_S$ takes its minimal value at $x\in S$, then $df(x)$
belongs to $N^*_0(S)$.

\medskip
(v) is nothing but Lemma 3.3 in \cite{KMS}.
\end{proof}

The following property of involutivity is due to J-M.~Bony
(\cite{B}).

\begin{theorem}\label{th:inv}
Let $S$ be a closed subset of $X$ and
$f$, $g$ two C$\,{}^1$-functions on $T^*X$ such that
$N^*_0(S)$ is contained in the zeroes's set of $f$ and that of $g$.
Then $N^*_0(S)$ is contained in the zeroes's set of
the Poisson bracket $\{f,g\}$.
\end{theorem}

\section{Involutivity of truncated microsupports}

\begin{theorem}[Weak involutivity of truncated microsupports]
Let $F\in D^b(\corps_X)$, $k\in\Z$ and let
$f$, $g$ be two C$\,{}^1$-functions on $T^*X$ such that
$\Ss_k(F)$ is contained in the zeroes's set of $f$ and that of $g$.
Then $\Ss_k(F)$ is contained in the zeroes's set of
the Poisson bracket $\{f,g\}$.
\end{theorem}

\begin{proof}
Assume that $\Ss_{k}(F)\subset \{f=0\}\cap\{g=0\}$ and let
$p=(x_0;\xi_0)$ such that $\{f,g\}(p)\not=0$.
We have to show that for any 
function $\phi$ with $\phi(x_0)=0$ and $d\phi(x_0)=\xi_0$,
the local cohomology $H^k_{\{\phi\ge0\}}(F)_{x_0}$ vanishes.

By induction on $k$
we may assume that $p\notin \Ss_{k-1}(F)$. 
Hence we may assume that $F\in D^{\ge k}(\corps_X)$.
Then 
$H^k_{\{\phi\ge0\}}(F)_{x_0}\simeq\Gamma_{\{\phi\ge0\}}(H^k(F))_{x_0}$.
Assume that $s\in \Gamma_{\{\phi\ge0\}}(H^k(F))_{x_0}$ does not vanish.
There exists an open neighborhood $V$ of $x_0$ such that
$s$ extends to $\tilde s\in \Gamma_{\{\phi\ge0\}}(V;H^k(F))$.
Then $S:=\supp(\tilde s)$ satisfies
$N^*_0(S)\subset \Ss_k(F)$. Hence $p\notin N^*_0(S)$
by Theorem \ref{th:inv}, which is a contradiction.
\end{proof}

\begin{remark}
The truncated microsupport is not involutive in the sense of
Definition \ref{def:invol}. Indeed,
for $X=\R$ and $Z =\{x\in X; x>0\}$, one has
\eqn
&&\Ss_0(\corps_{Z})=
\pset{(x;\xi)}{$\xi=0,\,x\ge0$}.
\eneqn
Hence one has $C_p(\Ss_0(F),\Ss_0(F))\subset\{-d\xi=0\}$ with $p=(0;0)$,
but $-\frac{\partial}{\partial x}=H(-d\xi)\notin C_p(\Ss_0(F))$.
\end{remark}

\begin{corollary}\label{co:setminus}
Let $S$ be a locally closed submanifold of $T^*X$
such that $T_pS$ is not involutive for any $p\in S$.
Then $S\cap\Ss_k(F)\subset \ol{\Ss_k(F)\setminus S}$.
\end{corollary}
\begin{theorem}\label{th:setminus}
Let $F\in D^b(\corps_X)$, $k\in\Z$ and let
$S$ be a subanalytic subset of $T^*X$
of dimension smaller than $\dim X$.
Then $\Ss_k(F)=\ol{\Ss_k(F)\setminus S}$.
\end{theorem}

\begin{proof}
We shall argue by the induction on $\dim S$.
There is a closed subanalytic subset $S_1$ of
$S$ such that $\dim S_1<\dim S_1$ and
$S_0:=S\setminus S_1$ is non singular.
Since $T_pS_0$ is not involutive for any $p\in S_0$,
$\Ss_k(F)\subset S_1\cup\ol{\Ss_k(F)\setminus S}$.
Hence one has $\ol{\Ss_k(F)\setminus S}=\ol{\Ss_k(F)\setminus S_1}$,
and the induction proceeds.
\end{proof}

\begin{theorem}\label{th:Rcons}
Let $F\in D^b_\Rc(\corps_X)$.
Let $\{Y_\alpha\}_{\alpha\in A}$ be a locally finite family of
real analytic submanifolds subanalytic in
$X$, and let $\Lambda_\alpha$ be an open subset of 
$T^*_{Y_\alpha}X$ subanalytic in $T^*X$,
such that
$\Ss(F)\subset\cup_{\alpha\in A}\ol{\Lambda_\alpha}$.
Let $K_\alpha\in D^b(\corps)$ and assume that
$F$ is microlocally isomorphic to
$\corps_{Y_\alpha}\otimes K_\alpha$
at every point of $\Lambda_\alpha$.
Set $A_k:=\{\alpha\in A;K_\alpha\notin D^{>k}(\corps)\}$.
Then for any $k\in\Z$,
$\Ss_k(F)=\bigcup\limits_{\alpha\in A_k}\ol{\Lambda_\alpha}$.
\end{theorem}
\begin{proof}
Set $S=\cup_{\alpha\in A}(\ol{\Lambda_\alpha}\setminus \Lambda_\alpha)$.
Then $\dim S<\dim X$.
If $\alpha\in A_k$, then
$\ol{\Lambda_\alpha}\subset\Ss_k(F)$
and if $\alpha\not\in A_k$, then $\Ss_k(F)\cap \Lambda_\alpha=\emptyset$.
Hence $\Ss_k(F)\setminus S=\Bigl(
\bigcup\limits_{\alpha\in A_k}\ol{\Lambda_\alpha}\Bigr)
\setminus S$.
Hence Theorem \ref{th:setminus} implies the desired result.
\end{proof}

The following corollary is proved in \cite{KMS}
when $\corps=\C$ by a different method.
Let $X$ be a complex manifold. 
Recall that $F\in D^b_\Cc(\corps_X)$ is perverse if
$$
\mbox{$\codim \Supp(H^k(F))\ge k$ and
$\codim \Supp(H^k(\rhom(F,\corps_X)))\ge k$}$$
for any $k\in\Z$.

\begin{corollary}
Let $X$ be a complex manifold. 
Let $F\in D^b_\Cc(\corps_X)$ and
let $\{X_\alpha\}_{\alpha\in A}$ be a family of
complex submanifolds such that
$\ol{X_\alpha}$ and $\ol{X_\alpha}\setminus X_\alpha$
are closed complex analytic subsets
and $\Ss(F)=\bigcup_{\alpha\in A} \ol{T^*_{X_\alpha}X}$.
\bnum
\item
If $F$ is a perverse sheaf,
then one has
\eq\label{eq:perv}
\Ss_k(F)=\bigcup_{\codim X_\alpha\le k}
\ol{T^*_{X_\alpha}X}\quad\mbox{for any $k$.}
\eneq
\item
Conversely if $F\in D^b(\corps_X)$
satisfies 
\eq\label{eq:pervs}\hspace{1cm}
\Ss_k(F)\cup\Ss_k(\rhom(F,\corps_X))
\subset\smash{\bigcup_{\codim X_\alpha\le k}}
\ol{T^*_{X_\alpha}X}\quad\mbox{for any $k$,}
\eneq
then $F$ is a perverse sheaf.
\enum
\end{corollary}
\begin{proof}
Recall (\cite[Theorem 10.3.12]{K-S1}) that $F$ is perverse if and only if 
$F$ is microlocally isomorphic to 
a finite direct sum of copies of $\corps_{X_\alpha}[-\codim X_\alpha]$
at a generic point of
$T^*_{X_\alpha}X$ for any $\alpha$. 

\noindent
(i) Assume $F$ is perverse. Using the notations of
Theorem \ref{th:Rcons}, we get 
\eqn
&&A_k=\{\alpha;\codim X_\alpha\le k\}
\eneqn
and the result follows from this theorem.

\noindent
(ii) The proof goes as \cite[Corollary 6.10]{KMS}. For the sake of
completeness, we repeat it. Recall that $\mu_Y(\cdot)$ denotes the
Sato's microlocalization functor.
By \cite{K-S1}, $F$ is isomorphic to 
$\corps_{X_\alpha}[-\codim X_\alpha]\otimes K$
at a generic point of
$T^*_{X_\alpha}X$ for some $K\in D^b(\corps)$.
Since $\mu_{X_\alpha}(F)$
must be in $D^{\ge \codim X_\alpha}(\corps_{T^*_{X_\alpha}X})$
and $\mu_{X_\alpha}(F)\simeq
\corps_{T^*_{X_\alpha}X}[-\codim X_\alpha]\otimes K$,
one has $K\in D^{\ge0}(\corps)$.
Similarly,
$\mu_{X_\alpha}\bl(\rhom(F,\corps_X)\br)\simeq
\corps_{T^*_{X_\alpha}X}[-\codim X_\alpha]\otimes \rhom(K,\corps)$
implies $K\in D^{\le0}(\corps)$.
\end{proof}

\medskip
{\small
\noindent
Masaki Kashiwara\\
Research Institute for Mathematical Sciences,
Kyoto University, Kyoto 606-8502 Japan\\
masaki@kurims.kyoto-u.ac.jp

\vspace{2mm}
\noindent
Teresa Monteiro Fernandes\\
Centro de {\'A}lgebra da Universidade de Lisboa, Complexo 2,\\
2 Avenida Prof. Gama Pinto, 1699 Lisboa codex Portugal\\
tmf@ptmat.lmc.fc.ul.pt

\vspace{2mm}
\noindent
Pierre Schapira\\
Universit{\'e} Pierre et Marie Curie, case 82,
Institut de Math{\'e}matiques,\\
4, place Jussieu, 75252 Paris cedex 05 France\\
schapira@math.jussieu.fr\\
http://www.institut.math.jussieu.fr/{\~{}}schapira/
}

\end{document}